\documentclass {article}
\usepackage{authblk}
\usepackage{etex}
\usepackage[utf8]{inputenc}
\usepackage[english]{babel}
\usepackage{amsmath}
\usepackage{amsthm}
\usepackage{amssymb}
\usepackage{sectsty}
\usepackage{titlesec}
\usepackage{color}
\usepackage[color,matrix,arrow]{xy}
\usepackage{amsgen}
\usepackage{amstext}
\usepackage{amsbsy}
\usepackage{amsopn}
\usepackage{amsfonts}
\usepackage{eepic}
\usepackage{graphicx}
\usepackage{epsf}
\usepackage{pstricks}
\usepackage{enumerate}
\usepackage{eqnarray}
\usepackage{faktor}
\xyoption{all}

\usepackage{pgf}
\usepackage{tikz-cd}
\usetikzlibrary{automata, arrows.meta, positioning,decorations.pathmorphing}

\newcommand{\footrecall}[1]{%
} 
\usetikzlibrary{snakes,shapes,arrows,automata}

\newcommand{\Mon}{\mbox{\rm Mon}}

\titleformat*{\section}{\large\bfseries}
\titleformat*{\subsection}{\normalsize \bfseries}

\newcommand{\N}{\mathbb{N}}
\newcommand{\Z}{\mathbb{Z}}

\newcommand{\Fix}{\text{Fix}}

\newcommand{\End}{\text{End}}
\newcommand{\Ker}{\text{Ker}}

\newcommand{\Eq}{\text{Eq}}
\newcommand{\Aut}{\text{Aut}}

\newcommand{\nin}{\not\in}

\theoremstyle{definition}
\newtheorem{theorem}{Theorem}[section]
\newtheorem{corollary}[theorem]{Corollary}

\newtheorem{proposition}[theorem]{Proposition}

\newtheorem{lemma}[theorem]{Lemma}

\begin{document}
 
 
\title{On the intersection of fixed subgroups of $F_n\times F_m$}
\author{Andr\'e Carvalho\\

Center for Mathematics and Applications (NOVA Math), NOVA SST

 2829–516 Caparica, Portugal

andrecruzcarvalho@gmail.com}
\maketitle

\begin{abstract}
We prove that, although it is undecidable if a subgroup fixed by an automorphism intersects nontrivially an arbitrary subgroup of $F_n\times F_m$, there is an algorithm that, taking as input a monomorphism and an endomorphism of $F_n\times F_m$, decides whether their fixed subgroups intersect nontrivially. The general case of this problem, where two arbitrary endomorphisms are given as input remains unknown. We show that, when two endomorphisms of a certain type are given as input, this problem is equivalent to the Post Correspondence Problem for free groups.
\end{abstract}

\section{Introduction}
Subgroups of elements fixed by endomorphisms have been an interesting object of study in group theory in the recent past. It started with the independent work of Gersten \cite{[Ger87]} and Cooper \cite{[Coo87]}, who proved that the subgroup of fixed points $\Fix(\varphi)$ of some automorphism $\varphi$ of $F_n$ is always finitely generated. Bestvina and Handel developed the theory of train tracks to prove that $\Fix(\varphi)$ has rank at most $n$ in \cite{[BH92]}. The problem of computing a basis for $\Fix(\varphi)$ was settled by Bogopolski and Maslakova in 2016 in \cite{[BM16]}. For general endomorphisms, Imrich and Turner \cite{[IT89]} introduced the concept of the stable image of an endomorphism, which is a useful tool to reduce problems on endomorphisms to problems on automorphisms, yielding in particular a proof generalizing the above mentioned results of Gersten, Cooper, Bestvina and Handel to general endomorphisms of the free group. Computability of the stable image (and so of the fixed subgroup) of an endomorphism was very recently proved by Mutanguha in \cite{[Mut22]}. Despite many positive results,  it is still not completely understood exactly which subgroups are fixed subgroups for some automorphism (resp. monomorphism, endomorphism) of a free group (see \cite{[Ven14]}). For example, it is not known whether we can decide if a given finitely generated subgroup of $F_n$ is fixed by some automorphism (or monomorphism or endomorphism) of $F_n$ \cite[Question 8.6]{[Ven14]}.

It is well known that arbitrary subgroups of the direct product of two free groups of finite rank can be complicated, and there are many undecidability results about them (see \cite{[BBN59],[Mil92], [Mih58]}). For example, by \cite[Theorem 4.5]{[Mil92]}, the isomorphism problem is undecidable for finitely generated subgroups of $F_2\times F_2$. Moreover, it is proved in \cite{[BR84]} that, even though there are continuously many non-isomorphic subgroups of $F_2\times F_2$, every finitely presented subgroup of a direct product of two free groups is a finite extension of a direct product of two free groups of finite rank.
 
We are interested in (a variation of) the \emph{Subgroup Intersection Triviality Problem (SITP($G$))}, which is the problem of deciding, on input two finite sets of generators for subgroups $H_1$, $H_2$ of $G$, whether $H_1\cap H_2$ is trivial or not.
As pointed out by Steinberg \cite{[Ste]}, Mihailova's construction \cite{[Mih58]}, together with a result in \cite{[BBN59]}, can be used to see that SITP($F_2\times F_2$) is undecidable. Indeed, 
given a finitely presented group $G=\langle X|R\rangle$ and a word $w$ in $X$, embedding $F_X\times F_X$ in $F_2\times F_2$, the intersection $$\{(u,v)\in F_X\times F_X \mid u=_Gv\}\cap \langle(1,w)\rangle$$ is trivial if and only if $w$ has infinite order in $G$, which is undecidable by \cite{[BBN59]}. So even the problem of deciding whether a finitely generated subgroup intersects a cyclic subgroup of $F_2\times F_2$ is undecidable. Throughout this paper, when we refer to a free group, we mean a free group of rank at least two. The case of free-abelian times free groups is easy since the intersection problem is solvable in this class of groups \cite{[DV13]}.

When a problem is not decidable for a class of groups, it might be possible to explore the  boundaries of its computability by considering special classes of subgroups.
The class of fixed subgroups of a group is, in general, a proper subclass of the class of all subgroups, but, as highlighted above  it is still unclear how \emph{special} they are.  Using Mutanguha's result for the free group and a classification of endomorphisms of $F_n\times F_m$ in seven different types, the author proved in \cite{[Car22b]} that there is an algorithm, taking an endomorphism $\varphi$ of $F_n\times F_m$ as input, to decide whether $\Fix(\varphi)$ is finitely generated or not and, in case it is, compute a finite set of generators for it. The proof gives a somewhat clear description of fixed subgroups giving hope of proving, for $F_n\times F_m$, decidability of the \emph{Fixed Subgroup Intersection Triviality Problem (FSITP($G$))}, a weaker version of the SITP($G$) that consists on deciding, on input two endomorphisms (resp. monomorphisms, automorphisms) $\varphi$, $\psi$ of $G$, whether $\Fix(\varphi)\cap \Fix(\psi)$ is trivial or not.

Fixed subgroups are, in principle, less complicated than arbitrary subgroups, so it might be possible to decide FSITP($F_n\times F_m$). However, we remark that, even when restricted to automorphisms, some complications can appear. Using the idea mentioned above, we can prove the following:
\newtheorem*{autogeneral}{Proposition \ref{autogeneral}}
\begin{autogeneral}
There is no algorithm that takes as input an automorphism $\Phi\in \Aut(F_n\times F_m)$ and a finitely generated subgroup $H\leq_{f.g.} F_n\times F_m$ and decides whether $\Fix(\Phi)\cap H$ is trivial or not.
\end{autogeneral}

Following this result, we prove that we can decide triviality of the intersection between the subgroup fixed by a monomorphism and the subgroup fixed by an endomorphism, showing that, in some sense, fixed subgroups are \emph{special} among all the other finitely generated subgroups.
\newtheorem*{cormain}{Corollary \ref{cormain}}
\begin{cormain}
There exists an algorithm with input natural numbers $m,n>1$, a monomorphism $\Phi\in \Mon(F_n\times F_m)$ and an endomorphism $\Psi\in \End(F_n\times F_m)$ that decides whether $\Fix(\Phi)\cap\Fix(\Psi)$ is trivial or not.
\end{cormain}

The Post Correspondence Problem (PCP) originally concerns free monoids and it is known to be undecidable \cite{[Pos46]}. Since its introduction, it has played an important role with many implications in computer science. The PCP has been studied for some classes of groups with positive results. For example, Myasnikov, Nikolaev and Ushakov proved in \cite[Theorem 5.8]{[MNU14]} that PCP($G$) belongs to \textbf{P} when $G$ is a finitely generated nilpotent group. For free groups, some results on the PCP and some of its variations have been obtained by Ciobanu and Logan (see \cite{[CL20],[CL21]}), but the general case remains unknown and an important problem (\cite[Problem 5.1.4]{[DKLM19]} \cite[Section 1.4]{[MNU14]}). We say that the \emph{PCP is decidable for free groups} if, given $m,n\in \N$ and homomorphisms $\varphi,\psi:F_m\to F_n$, we can decide if the equalizer $$\Eq(\varphi,\psi)=\{x\in F_m\mid x\varphi=x\psi\}$$
is trivial or not.

Corollary \ref{cormain} shows that fixed subgroups are easier than arbitrary subgroups, providing some hope of proving decidability of the general FSITP($F_n\times F_m$). However, we can prove that, despite being special, fixed subgroups of endomorphisms of $F_n\times F_m$ can still be complicated.  
\newtheorem*{posttypeIV}{Theorem \ref{posttypeIV}}
\begin{posttypeIV}
There exists an algorithm that takes as input two type IV endomorphisms $\Phi,\Psi\in\End(F_n\times F_m)$ and decides whether $\Fix(\Phi)\cap \Fix(\Psi)$ is trivial or not if and only if the PCP is decidable for free groups.
\end{posttypeIV}

This paper is organized as follows: in Section 2, we introduce the description  of endomorphisms of $F_n\times F_m$, following the work in  \cite{[Car22b]},  and some basic notation and results. In Section 3, we prove the main result: it is decidable whether the fixed subgroup of a monomorphism of $F_n\times F_m$ intersects the fixed subgroup of an endomorphism of $F_n\times F_m$. Finally, in Section 4, we discuss the general case where two endomorphisms (not necessarily injective) are considered, proving that for two endomorphisms of a certain type, this problem is equivalent to the PCP in free groups.

\section{Preliminaries}
The purpose of this section is to introduce the classification of endomorphisms of $F_n\times F_m$ obtained in  \cite{[Car22b]}, as well two important results on rational subsets of groups.

 For $(i,j)\in[n]\times [m]$, we define   $\lambda_i:F_n\to \mathbb Z$ as the endomorphism given by $a_k\mapsto \delta_{ik}$ and $\tau_j:F_m\to \mathbb Z$ given by $b_k\mapsto \delta_{jk}$, where $\delta_{ij}$ is the Kronecker symbol.
 
For $x\in F_n, y\in F_m$ and integers $p_i,q_i,r_j,s_j\in \Z$,we denote $\sum\limits_{i\in[n]} \lambda_i(x)p_i$ by $x^{P}$; $\sum\limits_{j\in[m]} \tau_j(y)r_j$ by $y^{R}$; $\sum\limits_{i\in[n]} \lambda_i(x)q_i$ by $x^{Q}$ and $\sum\limits_{j\in[n]} \tau_j(y)s_i$ by $y^{S}$. 
We also define $P=\{p_i\in\Z\mid i\in[n]\}$, $Q=\{q_i\in\Z\mid i\in[n]\}$, $R=\{r_j\in\Z\mid j\in[m]\}$ and $S=\{s_j\in\Z\mid j\in[m]\}$.
 We will keep this notation throughout the paper.

In \cite{[Car22b]}, the author classified endomorphisms of the direct product of two finitely generated free groups $F_n\times F_m$, with $m,n>1$ in seven different types:
\begin{enumerate}[(I)]
\item $(x,y)\mapsto\left(u^{x^P+y^R},v^{x^Q+y^S}\right)$, for some $1\neq u \in F_n$, $1\neq v \in F_m$ and integers $p_i,q_i, r_j, s_j\in \mathbb Z$ for $(i,j)\in [n]\times [m]$, such that $P,Q,R,S\neq \{0\}$. 
\item $(x,y)\mapsto\left(y\phi,v^{x^Q+y^S}\right)$, for some nontrivial homomorphism $\phi:F_m\to F_n$, $1\neq v \in F_m$  and integers $q_i,  s_j\in \mathbb Z$ for $(i,j)\in [n]\times [m]$,  such that $Q,S\neq \{0\}$. 
\item $(x,y)\mapsto\left(u^{x^P+y^R},y\phi\right),$  for some nontrivial endomorphism $\phi\in\End(F_m)$, $1\neq u \in F_n$, and integers $p_i,  r_j\in \mathbb Z$ for $(i,j)\in [n]\times [m]$,  such that $P,R\neq \{0\}$. 
\item $(x,y)\mapsto(y\phi,y\psi)$, for some nontrivial homomorphism $\phi: F_m\to F_n$ and nontrivial endomorphism $\psi\in \End(F_m)$.
\item $(x,y)\mapsto\left(1,v^{x^Q+y^S}\right)$, for some $1\neq v \in F_m$, and integers $q_i,  s_j\in \mathbb Z$ for $(i,j)\in [n]\times [m]$, such that $Q,S\neq \{0\}$.
\item $(x,y)\mapsto (x\phi,y\psi)$, for some endomorphisms $\phi\in \End(F_n)$, $\psi\in \End(F_m)$.
\item $(x,y)\mapsto (y\psi,x\phi)$, for homomorphisms $\phi:F_n\to F_m$  and $\psi:F_m\to F_n$.
\end{enumerate}

 From \cite[Proposition 3.2]{[Car22b]}, injective endomorphisms correspond to endomorphisms of type VI or VII such that the component mappings $\phi$ and $\psi$ are injective.  

Additionally, the author studied fixed subgroups for endomorphisms of each of the seven types, obtaining structural results and an algorithm to decide whether the fixed subgroup of an endomorphism is finitely generated or not and compute it in case it is. In doing so, it was convenient to split the type III endomorphisms in two: we say that $\Phi\in \End(F_n\times F_m)$ is of type III.1 if it is of type III and  $u^P\neq 1$ and of type III.2 if it is of type III but  $u^P= 1$.

Now we introduce the notion of a rational subset of a group and two important results concerning rational subsets.

Given a finitely generated group $G=\langle A\rangle$, a finite generating set $A$ and a set of formal inverses $A^{-1}$, write $\tilde A=A\cup A^{-1}$. There is a canonical (surjective) homomorphism $\pi:\tilde A^*\to G$
mapping $a\in \tilde A$ (resp. $a^{-1}\in \tilde A$) to  $a\in  G$ (resp. $a^{-1}\in  G$).

A subset $K\subseteq G$ is said to be \emph{rational} if there is some rational language $L\subseteq \tilde A^*$ such that $L\pi=K$.
The set of all such subsets is denoted by $\text{Rat } G$. Rational subsets generalize the notion of finitely generated subgroups.

\begin{theorem}\cite[Theorem III.2.7]{[Ber79]}
\label{anisimovseifert}
Let $H$ be a subgroup of a group $G$. Then $H\in \text{Rat } G$ if and only if $H$ is finitely generated.
\end{theorem}

In case the group $G$ is a free group with basis $A$, we define the set of reduced words  of $L\subseteq \tilde A^*$ by
$$\overline L=\{\bar w \mid w\in L \}.$$
Benois' Theorem provides us with a useful characterization of rational subsets of free groups in terms of reduced words representing the elements in the subset.

\begin{theorem}\cite[Benois' theorem]{[Ben79]}
\label{benois}
Let $F$ be a free group with basis $A$. Then, a subset of $\overline{\tilde A^*}$ is a rational language of $\tilde A^*$ if and only if it is a rational subset of $F$.
\end{theorem}

\section{The main result} 
Using the same idea from \cite{[Ste]} and the description of automorphisms of $F_n\times F_m$ from \cite{[Car22b]}, we can prove the following:
\begin{proposition}
\label{autogeneral}
There is no algorithm that takes as input an automorphism $\Phi\in \Aut(F_n\times F_m)$ and a finitely generated subgroup $H\leq_{f.g.} F_n\times F_m$ and decides whether $\Fix(\Phi)\cap H$ is trivial or not.
\end{proposition}
\noindent\textit{Proof.} Let $G=\langle X|R\rangle$ be a finitely presented group and  $w$ be a word in $X$. Write $w$ as $z^k$, where $z$ is a word in $X$ that is not a proper power (if $w$ is not a proper power itself, then we take $z=w$). Consider the automorphism $\Phi$ of $F_X\times F_X$ defined by $(x,y)\mapsto (x\varphi,y\lambda_z)$, where $\varphi$ is a letter permutation automorphism with no fixed letters and $\lambda_z$ is the inner automorphism of $F_X$ given by conjugation by $z$. Then $\Fix(\lambda_z)=C_{F_X}(z)=\langle z \rangle$ and $$\Fix(\Phi)=\Fix(\varphi)\times\Fix(\lambda_z)=\{(1,z^k)\mid k\in \Z\}=\langle(1,z)\rangle.$$

 We have that $w$ has infinite order in $G$ if and only if $z$ has infinite order in $G$. The subgroup $\{(u,v)\in F_X\times F_X \mid u=_Gv\}$ is finitely generated and the
 intersection $$\{(u,v)\in F_X\times F_X \mid u=_Gv\}\cap \langle(1,z)\rangle=\{(u,v)\in F_X\times F_X \mid u=_Gv\}\cap \Fix(\Phi)$$ is trivial if and only if $z$ has infinite order in $G$, which is undecidable by \cite{[BBN59]}. \qed\\

Now, we present a technical lemma that will be used often in the proof of the main result.
\begin{lemma}\label{lem: power fixed}
Let $\varphi\in \End(F_n)$ and $u\in F_n\setminus \{1\}$, $a\in \Z\setminus \{0\}$ such that $u^a\in \Fix(\varphi)$. Then $u\in \Fix(\varphi)$.  
\end{lemma}
\noindent\textit{Proof.} We have that $(u\varphi)^a=u^a\varphi=u^a$. Writing $u\varphi=wzw^{-1}$ and $u=xyx^{-1}$, where $z$ and $y$ are the cyclically reduced cores of $u\varphi$ and $u$, respectively, we have that $(u\varphi)^a=wz^aw^{-1}$ and $u^a=xy^ax^{-1}$ and so $w=x$ and $z^a=y^a$ which means that $z=y$. So $u=u\varphi$. 
\qed\\

We are not able to solve the general case, when two arbitrary endomorphisms are given as input. However, we prove decidability in the case where one of the endomorphisms is of type VI or VII, which includes in particular all injective endomorphisms. This shows that fixed subgroups are algorithmically easier than arbitrary subgroups.

\begin{theorem}
There exists an algorithm with input natural numbers $m,n>1$, an endomorphism $\Phi\in \End(F_n\times F_m)$ of type VI or VII and an endomorphism $\Psi\in \End(F_n\times F_m)$ that decides whether $\Fix(\Phi)\cap\Fix(\Psi)$ is trivial or not.
\end{theorem}
\noindent\textit{Proof.} 
We start by considering the case where $\Phi$ is a type VI endomorphism and then the case where $\Phi$ has type VII.

\textbf{Case 1: $\Phi$ is of type VI.}
We have that $\Phi$ is given by  $(x,y)\mapsto (x\phi,y\psi)$, for some $\phi\in\End(F_n)$ and $\psi\in \End(F_m)$, and so, $\Fix(\Phi)=\Fix(\phi)\times \Fix(\psi)$, 
We now consider several subcases. 

\textbf{Subcase 1.1: $\Psi$ is of type I} We have that $\Psi$ is defined by 
$$(x,y)\mapsto\left(u^{x^P+y^R},v^{x^Q+y^S}\right).$$ As proved in \cite[Subsection 4.1]{[Car22b]},
we have that  $\Fix(\Psi)=\{(u^a,v^b)\in F_n\times F_m \mid (a,b)\in \Ker(M_\Psi)\}$, where  $$M_\Psi=\begin{bmatrix} -1+u^P & v^R\\
u^Q& -1+v^S
 \end{bmatrix}.$$
So, we check whether $u\in \Fix(\phi)$ and $v\in \Fix(\psi)$. In view Lemma \ref{lem: power fixed}, if $u\not\in \Fix(\phi)$ (resp. $v\not\in \Fix(\psi)$), then $u^a\in \Fix(\phi)$ (resp. $v^b\in \Fix(\psi)$) if and only if $a=0$ (resp. $b=0$).
 If $u\in\Fix(\phi)$ and $v\in\Fix(\psi)$, then $\Fix(\Phi)\cap \Fix(\Psi)=\Fix(\Psi)$, and it is computable. If $u\not\in\Fix(\phi)$ and $v\not\in\Fix(\psi)$, then $\Fix(\Phi)\cap \Fix(\Psi)$ is trivial. If, say $u\nin\Fix(\phi)$ but $v\in\Fix(\psi)$, then 
 \begin{align*}
 \Fix(\Phi)\cap \Fix(\Psi)&=\{(1,v^b)\in F_n\times F_m \mid (0,b)\in \Ker(M_\Psi)\}\\&=\{(1,v^b)\in F_n\times F_m \mid bv^R=-b+bv^S=0\},
 \end{align*}
 whose triviality can be decided (it is trivial if $v^S\neq 1$ or $v^R\neq 0$ and infinite cyclic otherwise). The case where $u\in\Fix(\phi)$ but $v\nin\Fix(\psi)$ is analogous.

\textbf{Subcase 1.2: $\Psi$ is of type II.} We have that $\Psi$ is defined by 
$$(x,y)\mapsto\left(y\theta,v^{x^Q+y^S}\right),$$
where $\theta:F_m\to F_n$ is a homomorphism.
 As proved in \cite[Subsection 4.2]{[Car22b]}, if $(v\theta)^Q +v^S\neq 1$, we have that $\Fix(\Psi)$ is trivial (and so is $\Fix(\Phi)\cap \Fix(\Psi)$). If not, $\Fix(\Psi)$ is given by $\{(v^b\theta,v^b)\in F_n\times F_m \mid b\in \mathbb Z\}$. In this case, to decide whether 
 $\Fix(\Phi)\cap \Fix(\Psi)$ is trivial or not, we start by checking if $v\in \Fix(\psi)$. If not, then the intersection is trivial since $v^b\in \Fix(\psi)$ if and only if $b=0$. If $v\in \Fix(\psi)$, then
  \begin{align*}
 \Fix(\Phi)\cap \Fix(\Psi)&=\{(v^b\theta,v^b)\in F_n\times F_m \mid v^b\theta\in \Fix(\phi)\}\\
 &\simeq\{b\in \Z \mid v^b\theta\in \Fix(\phi)\}\\
&\simeq \langle v\theta\rangle\cap \Fix(\phi).
 \end{align*}
 Since we can compute $ \langle v\theta\rangle\cap \Fix(\phi)$, we can decide its triviality.
 
\textbf{Subcase 1.3: $\Psi$ is of type III.1.} We have that $\Psi$ is defined by 
$$(x,y)\mapsto\left(u^{x^P+y^R},y\theta\right),$$ with $u^P\neq 1$, and, as shown in \cite[Subsection 4.3]{[Car22b]},
 putting $$H=\left\{y\in \Fix(\theta)\;\bigg\lvert\; \left(1-u^P\right) \text{ divides } {y^R}\right\},$$ we have that
$$\Fix(\Psi)=\left\{\left(u^{\frac{y^R}{1- u^P}},y\right)\mid y\in H \right\}.$$
Also, in \cite[Subsection 4.3]{[Car22b]} it is proved that $H$ is the intersection of $\Fix(\theta)$ with a finite index subgroup of $F_m$ and can be computed. So, we start by checking if $H\cap \Fix(\psi)$ is trivial or not. If it is, then the only possible value for $y$ is $1$, which makes $u^{\frac{y^R}{1- u^P}}$ trivial. Thus, in this case $\Fix(\Phi)\cap\Fix(\Psi)=1$. Assume then that $H\cap \Fix(\psi)$ is nontrivial. If $u\in \Fix(\phi)$, then
 $$\Fix(\Phi)\cap\Fix(\Psi)=\left\{\left(u^{\frac{y^R}{1- u^P}},y\right)\mid y\in H\cap \Fix(\psi) \right\}\simeq H\cap\Fix(\psi)\neq 1.$$
 If, on the other hand, $u\not\in \Fix(\phi)$, then
 $$\Fix(\Phi)\cap\Fix(\Psi)=\left\{\left(1,y\right)\mid y\in H\cap \Fix(\psi) : y^R=0\right\}.$$
 In  \cite[Subsection 4.3]{[Car22b]}, it is proved that the language of reduced words $y$ in $F_m$ such that $y^R=0$ is context-free. Also, by Theorems \ref{anisimovseifert} and \ref{benois}, the language of reduced words representing elements in $H\cap \Fix(\psi)$ is rational. So the intersection of both languages is an effectively constructible context-free language. So, $\Fix(\Phi)\cap\Fix(\Psi)$ is isomorphic to an (effectively constructible) context-free subgroup of $F_m$ and its triviality can be decided, since it corresponds to deciding whether this context-free language contains only the empty word.

\textbf{Subcase 1.4: $\Psi$ is of type III.2.} This case is similar to the previous. We have that $\Psi$ is defined by 
$$(x,y)\mapsto\left(u^{x^P+y^R},y\theta\right),$$ with $u^P= 1$, and, as shown in \cite[Subsection 4.3]{[Car22b]},
 putting $H=\left\{y\in \Fix(\theta)\mid y^R=0\right\}$, we have that 
$\Fix(\Psi)=\{(u^a,y)\mid y\in H, a\in\Z \}$. Again, if $u\not\in\Fix(\phi)$, then 
$$\Fix(\Phi)\cap\Fix(\Psi)=\{(1,y)\mid y\in H\cap \Fix(\psi)\},$$
which is isomorphic to a context-free subgroup of $F_m$. If $u\in \Fix(\phi)$, then
$$\Fix(\Phi)\cap\Fix(\Psi)=\{(u^a,y)\mid a\in \Z, \, y\in H\cap \Fix(\psi)\},$$
which is isomorphic to $\Z\times (H\cap \Fix(\psi))$, the direct product of $\Z$ with a context-free subgroup of $F_m$. In particular, in this case, $\Fix(\Phi)\cap\Fix(\Psi)$ is nontrivial.

\textbf{Subcase 1.5: $\Psi$ is of type IV.}  We have that $\Psi$ is given by $(x,y)\mapsto(y\theta,y\sigma),$
where $\theta:F_m\to F_n$ is a homomorphism and $\sigma\in \End(F_m)$. As shown in \cite[Subsection 4.4]{[Car22b]},
$\Fix(\Psi)=\{(y\theta,y)\in F_n\times F_m \mid y\in\Fix(\sigma)\}$. Hence, 
$$\Fix(\Phi)\cap\Fix(\Psi)=\{(y\theta,y)\in F_n\times F_m \mid y\in\Fix(\sigma)\cap\Fix(\psi)\cap\Fix(\phi)\theta^{-1}\}.$$
We compute $(\Fix(\sigma)\cap\Fix(\psi))\theta \cap \Fix(\phi).$ If it is nontrivial, then there is some $y\in \Fix(\sigma)\cap\Fix(\psi)$ such that $y\theta\in \Fix(\phi)$ is nontrivial, so $(y\theta,y)$ is a nontrivial element of $\Fix(\Phi)\cap\Fix(\Psi)$. If it is trivial, then
$$\Fix(\Phi)\cap\Fix(\Psi)=\{(1,y)\in F_n\times F_m \mid y\in\Fix(\sigma)\cap\Fix(\psi)\cap\Ker(\theta)\},$$
and it is trivial if and only if $\theta|_{\Fix(\sigma)\cap\Fix(\psi)}$ is injective. This can be checked: $\Fix(\sigma)\cap\Fix(\psi)$ is a computable finitely generated subgroup of $F_m$ and so $(\Fix(\sigma)\cap\Fix(\psi))\theta$ is a computable finitely generated subgroup of $F_n$, hence they are both free groups of finite rank, and their ranks can be computed. If they do not coincide, then $\theta|_{\Fix(\sigma)\cap\Fix(\psi)}$ is not injective, since the two subgroups are not isomorphic; if they coincide, then 
$$\Fix(\sigma)\cap\Fix(\psi)\simeq(\Fix(\sigma)\cap\Fix(\psi))\theta\simeq\faktor{\Fix(\sigma)\cap\Fix(\psi)}{\Ker(\theta|_{\Fix(\sigma)\cap\Fix(\psi)})}$$ and so, since free groups of finite rank are hopfian, $\Ker(\theta|_{\Fix(\sigma)\cap\Fix(\psi)})$ must be trivial.

\textbf{Subcase 1.6: $\Psi$ is of type V.} We have that $\Psi$ is defined by 
$$(x,y)\mapsto\left(1,v^{x^Q+y^S}\right).$$
It is shown in \cite[Subsection 4.5]{[Car22b]} that, if $v^S\neq 1$, then $\Fix(\Psi)$ is trivial (and so is $\Fix(\Phi)\cap\Fix(\Psi)$) and if $v^S=1$, then $\Fix(\Psi)=\{(1,v^b)\mid b\in\Z\}$. Hence, if $v\in \Fix(\psi)$, then
$$\Fix(\Phi)\cap\Fix(\Psi)=\{(1,v^b)\mid b\in\Z\}$$ and if not, then $\Fix(\Phi)\cap\Fix(\Psi)$ is trivial.

\textbf{Subcase 1.7: $\Psi$ is of type VI.} In this case, $\Psi$ is given by  $(x,y)\mapsto (x\theta,y\sigma)$, for some $\theta\in\End(F_n)$ and $\sigma\in \End(F_m)$. Then it is clear that $\Fix(\Psi)=\Fix(\theta)\times \Fix(\sigma)$ and 
$$\Fix(\Phi)\cap \Fix(\Psi)=(\Fix(\phi)\cap\Fix(\theta))\times (\Fix(\psi)\cap\Fix(\sigma)),$$
which can be effectively computed.

\textbf{Subcase 1.8: $\Psi$ is of type VII.} This case is similar to Subcase 1.5. We have that $\Psi$ is given by $(x,y)\mapsto (y\sigma,x\theta)$, for homomorphisms $\theta:F_n\to F_m$  and $\sigma:F_m\to F_n$. As proved in \cite[Subsection 4.7]{[Car22b]},
$$\Fix(\Psi)=\{(x,x\theta)\mid x\in \Fix(\theta\sigma)\}=\{(y\sigma,y)\mid y\in \Fix(\sigma\theta)\}.$$
So,
$$\Fix(\Phi)\cap\Fix(\Psi)=\{(x,x\theta)\mid x\in \Fix(\theta\sigma)\cap\Fix(\phi)\cap(\Fix(\psi))\theta^{-1}\}.$$
Proceeding as in the Subcase 1.5, we can check its triviality.

\textbf{Case 2: $\Phi$ is of type VII.} We have that $\Phi$ is given by $(x,y)\mapsto (y\psi,x\phi)$, for homomorphisms $\phi:F_n\to F_m$  and $\psi:F_m\to F_n$ and, as mentioned above, 
$$\Fix(\Phi)=\{(x,x\phi)\mid x\in \Fix(\phi\psi)\}=\{(y\psi,y)\mid y\in \Fix(\psi\phi)\}.$$

\textbf{Subcase 2.1: $\Psi$ is of type I.} Recall that $$\Fix(\Psi)=\{(u^a,v^b)\in F_n\times F_m \mid (a,b)\in \Ker(M_\Psi)\}.$$ We check if $u\in \Fix(\phi\psi).$ If  not, then the only possible value for $x$ is $1$, and 
$\Fix(\Phi)\cap \Fix(\Psi)$ must be trivial. If, on the other hand,  $u\in \Fix(\phi\psi)$, then 
$$\Fix(\Phi)\cap \Fix(\Psi)=\{(u^a, v^b)\in F_n\times F_m\mid (a,b)\in \Ker(M_\Psi) \text{ and } v^b=u^a\phi\}$$
We compute the solutions $(m,n)$ such that $v^m=(u\phi)^n$, which is a subgroup of $\Z^2$ of the form $\langle (m,n)\rangle$, for computable $m,n\in \Z$. Then $\Fix(\Phi)\cap \Fix(\Psi)$ is trivial if and only if $\Ker(M_\Psi)\cap \langle (m,n)\rangle$ is trivial.

\textbf{Subcase 2.2: $\Psi$ is of type II.} We have that $\Psi$ is defined by 
$$(x,y)\mapsto\left(y\theta,v^{x^Q+y^S}\right),$$
 and if $(v\theta)^Q +v^S\neq 1$, we have that $\Fix(\Psi)$ is trivial (and so is $\Fix(\Phi)\cap \Fix(\Psi)$). If not, $\Fix(\Psi)$ is given by $\{(v^b\theta,v^b)\in F_n\times F_m \mid b\in \mathbb Z\}$. In this case, to decide whether 
 $\Fix(\Phi)\cap \Fix(\Psi)$ is trivial or not, we start by checking if $v\in \Fix(\psi\phi)$. If not, then  $\Fix(\Phi)\cap \Fix(\Psi)$ must be trivial. If  $v\in \Fix(\psi\phi)$, then 
\begin{align*}
\Fix(\Phi)\cap \Fix(\Psi)&=\{(v^b\theta,v^b)\mid b\in\Z \text{ such that } v^b\theta =v^b\psi\}\\
&=
\begin{cases}
\{(v^b\theta,v^b)\mid b\in\Z \} \quad &\text{if $v\theta =v\psi$}\\
1 \quad &\text{if $v\theta \neq v\psi$}
\end{cases}.
\end{align*}

\textbf{Subcase 2.3: $\Psi$ is of type III.1.}  We have that $\Psi$ is defined by 
$(x,y)\mapsto\left(u^{x^P+y^R},y\theta\right),$ with $u^P\neq 1$, and, 
$$\Fix(\Psi)=\left\{\left(u^{\frac{y^R}{1- u^P}},y\right)\mid y\in \Fix(\theta) \text{ such that }  \left(1-u^P\right) \text{ divides } {y^R}  \right\}.$$
We start by checking whether $u\in \Fix(\phi\psi)$ or not.  If not, then $\Fix(\Phi)\cap\Fix(\Psi)$ must be trivial. If $u\in \Fix(\phi\psi)$, then we check if $u\phi\in \Fix(\theta)$. If not, then $\Fix(\Phi)\cap\Fix(\Psi)$ must be trivial, since the first component must be a power of $u$ and the second must be a power of $u\phi$ and $(u\phi)^k\in \Fix(\theta)$ implies that $k=0$. So, in this case an element of the intersection must have trivial second component and so it must be trivial. We now have to deal with the case where $u\in \Fix(\phi\psi)$ and  $u\phi\in \Fix(\theta)$. In this case, the candidates are of the form $(u^k, u^k\phi)$ and all of these elements belong to $\Fix(\Phi)$.  Since $(u^k\phi)^R=k(u\phi)^R$, the elements in the intersection are the points $(u^k, u^k\phi)$ for $k\in \Z$ such that $$k=\frac{k(u\phi)^R}{1-u^P}.$$
So, if $\frac{(u\phi)^R}{1-u^P}\neq 1$, then $\Fix(\Phi)\cap\Fix(\Psi)$ is trivial and if, on the other hand, $\frac{(u\phi)^R}{1-u^P}= 1$, then  $$\Fix(\Phi)\cap\Fix(\Psi)=\{(u^k,u^k\phi)\in F_n\times F_m \mid k\in \Z\}.$$

\textbf{Subcase 2.4: $\Psi$ is of type III.2} We have that $\Psi$ is defined by 
$$(x,y)\mapsto\left(u^{x^P+y^R},y\theta\right),$$ with $u^P= 1$, and
$\Fix(\Psi)=\{(u^a,y)\mid y\in  \Fix(\theta),\, y^R=0,\, a\in\Z \}$. We check if $u\in \Fix(\phi\psi)$. If not, then $\Fix(\Phi)\cap \Fix(\Psi)$ is trivial. If $u\in \Fix(\phi\psi)$, then $\Fix(\Phi)\cap\Fix(\Psi)$ consists of the elements of the form $(u^a,u^a\phi)$ such that $u^a\phi\in \Fix(\theta)$ and $(u^a\phi)^R=a(u\phi)^R=0$. We compute $(u\phi)^R$ and if $(u\phi)^R\neq 0$, then $a$ must be $0$ and the intersection must be trivial. If $(u\phi)^R= 0$, then if $u\phi\not\in \Fix(\theta)$, $\Fix(\Phi)\cap\Fix(\Psi)$ is trivial and otherwise 
$\Fix(\Phi)\cap \Fix(\Psi)=\{(u^a,u^a\phi)\mid  a\in\Z \}.$

\textbf{Subcase 2.5: $\Psi$ is of type IV.} We have that $\Psi$ is given by $(x,y)\mapsto(y\theta,y\sigma),$
where $\theta:F_m\to F_n$ is a homomorphism and $\sigma\in \End(F_m)$ and
$\Fix(\Psi)=\{(y\theta,y)\in F_n\times F_m \mid y\in\Fix(\sigma)\}$.

 Hence, 
$$\Fix(\Phi)\cap\Fix(\Psi)=\{(x,x\phi)\in F_n\times F_m \mid x\in\Fix(\phi\psi)\cap\Fix(\phi\theta)\cap\Fix(\sigma)\phi^{-1}\}.$$

 We start by computing $(\Fix(\phi\psi)\cap\Fix(\phi\theta))\phi \cap \Fix(\sigma).$ If it is nontrivial, then there is some $y\in \Fix(\phi\psi)\cap\Fix(\phi\theta)$ such that $y\phi \in \Fix(\sigma)$ is nontrivial, so $(y,y\phi)$ is a nontrivial element of $\Fix(\Phi)\cap\Fix(\Psi)$.
 If it is trivial, then the second component of an element in $\Fix(\Phi)\cap \Fix(\Psi)$ must be $1$, since $1$ is the only fixed point of $\Fix(\sigma)$, which is also the image by $\phi$ of a point fixed by both $\phi\psi$ and $\phi\theta$ (and the projection to the first component must be contained in $\Fix(\phi\psi)\cap \Fix(\phi\theta)$). But any element in $\Fix(\Psi)$ with trivial second component must also have trivial first component. So $\Fix(\Phi)\cap\Fix(\Psi)$ is trivial if and only if  $(\Fix(\phi\psi)\cap\Fix(\phi\theta))\phi \cap \Fix(\sigma)$ is trivial and that is decidable.

\textbf{Subcase 2.6: $\Psi$ is of type V.}  We have that $\Psi$ is defined by 
$(x,y)\mapsto\left(1,v^{x^Q+y^S}\right)$ and,  if $v^S\neq 1$, then $\Fix(\Psi)$ is trivial (and so is $\Fix(\Phi)\cap\Fix(\Psi)$);  if $v^S=1$, then $\Fix(\Psi)=\{(1,v^b)\mid b\in\Z\}$. Hence,  if $v\nin \Fix(\psi\phi)$, then $\Fix(\Phi)\cap\Fix(\Psi)$ is trivial and if $v\in \Fix(\psi\phi)$, then
$$\Fix(\Phi)\cap\Fix(\Psi)=\{(1,v^b)\mid b\in\Z \text{ such that } v^b\in\Ker(\psi)\}$$ and 
so it is trivial if and only if $v\in \Ker(\psi)$ and equal to $\{(1,v^b)\mid b\in\Z\}$
otherwise.

\textbf{Subcase 2.7: $\Psi$ is of type VI.} This is already proved in Subcase 1.8.

\textbf{Subcase 2.8: $\Psi$ is of type VII.} We want to check whether 
\begin{align*}
\Fix(\Phi)\cap \Fix(\Psi)&=\{(x,x\phi)\mid x\in\Fix(\phi\psi)\}\cap \{(y\sigma,y)\mid y\in\Fix(\sigma\theta)\}\\
&= \{(y\sigma,y)\mid y\in\Fix(\sigma\theta)\cap \Fix(\sigma\phi)\cap \Fix(\phi\psi)\sigma^{-1}\}
\end{align*}
 is trivial or not. This case is similar to Subcase 2.5: we compute $(\Fix(\sigma\theta)\cap\Fix(\sigma\phi))\sigma \cap \Fix(\phi\psi)$, which is trivial if and only if $\Fix(\Phi)\cap \Fix(\Psi)$ is trivial too. 
\qed

\begin{corollary}\label{cormain}
There exists an algorithm with input natural numbers $m,n>1$, a monomorphism $\Phi\in \Mon(F_n\times F_m)$ and an endomorphism $\Psi\in \End(F_n\times F_m)$ that decides whether $\Fix(\Phi)\cap\Fix(\Psi)$ is trivial or not.
\end{corollary}

\section{The general case}

We are not able to solve this problem in its complete generality considering general endomorphisms. We now show that this is a difficult problem, by showing that, when two endomorphisms of type IV are considered, this problem is equivalent to the PCP for free groups. 

\begin{theorem}\label{posttypeIV}
There exists an algorithm that takes as input two type IV endomorphisms $\Phi,\Psi\in\End(F_n\times F_m)$ and decides whether $\Fix(\Phi)\cap \Fix(\Psi)$ is trivial or not if and only if the PCP is decidable for free groups.
\end{theorem}
\noindent\textit{Proof.} Suppose that the subgroup intersection triviality problem is decidable for fixed subgroups of type IV endomorphisms of $F_n\times F_m$ and let $\phi,\psi:F_n\to F_m$ be two homomorphisms. Define two type IV endomorphisms $\Phi,\Psi\in\End(F_n\times F_m)$ by $(x,y)\mapsto (y\phi,y)$ and $(x,y)\mapsto (y\psi,y)$, respectively. Then $\Fix(\Phi)=\{(y\phi,y)\in F_n\times F_m \mid y\in F_m\}$ and $\Fix(\Psi)=\{(y\psi,y)\in F_n\times F_m \mid y\in F_m\}$. Hence,
$$\Fix(\Phi)\cap \Fix(\Psi)=\{(y\phi,y)\in F_n\times F_m \mid y\phi=y\psi\}\simeq \Eq(\phi,\psi),$$
and, since we can decide if $\Fix(\Phi)\cap \Fix(\Psi)$ is trivial, we can decide if $\Eq(\phi,\psi)$ is trivial.

Conversely, suppose that the PCP is decidable for free groups. Let $\Phi,\Psi\in \End(F_n\times F_m)$ be two type IV endomorphisms of $F_n\times F_m$. Then, there are homomorphisms $\phi,\theta:F_m\to F_n$ and endomorphisms $\psi,\sigma\in \End(F_m)$ such that $\Phi$ and $\Psi$ are given by $(x,y)\mapsto(y\phi,y\psi)$ and $(x,y)\mapsto(y\theta,y\sigma),$ respectively. Moreover, $\Fix(\Phi)=\{(y\phi,y)\in F_n\times F_m \mid y\in\Fix(\psi)\}$ and $\Fix(\Psi)=\{(y\theta,y)\in F_n\times F_m \mid y\in\Fix(\sigma)\}$, so
\begin{align*}
\Fix(\Phi)\cap\Fix(\Psi)&=\{(y\theta,y)\in F_n\times F_m \mid y\in\Fix(\psi)\cap\Fix(\sigma)\cap \Eq(\phi,\theta)\}\\
&\simeq \Fix(\psi)\cap\Fix(\sigma)\cap \Eq(\phi,\theta).
\end{align*}
We can decide if $\Fix(\psi)\cap\Fix(\sigma)\cap \Eq(\phi,\theta)$ is trivial using the PCP since ${\Fix(\psi)\cap\Fix(\sigma)}$ is a finitely generated (and computable) free group and
 $$\Fix(\psi)\cap\Fix(\sigma)\cap \Eq(\phi,\theta)=\Eq(\phi|_{\Fix(\psi)\cap\Fix(\sigma)},\theta|_{\Fix(\psi)\cap\Fix(\sigma)}).$$ 
\qed\\
\section*{Acknowledgements}
The author is grateful to Enric Ventura for suggesting the study of FSITP($F_n\times F_m$). This work is funded by national funds through the FCT - Fundação para a Ciência e a Tecnologia, I.P., under the scope of the projects UIDB/00297/2020 and UIDP/00297/2020 (Center for Mathematics and Applications).

\bibliographystyle{plain}
\bibliography{Bibliografia}

 \end{document}